\documentclass[12pt]{article}
\usepackage{graphpap}
\usepackage{stmaryrd}
\usepackage{amssymb}
\usepackage{amsmath}
\usepackage{algorithm}
\usepackage[noend]{algorithmic}
\usepackage{bm}

\newtheorem{lem}{Lemma}
\newtheorem{Theorem}{Theorem}
\newtheorem{defi}{Definition}
\newtheorem{crl}{Corollary}
\newtheorem{prp}{Proposition}
\newtheorem{rmk}{Remark}
\newtheorem{xmp}{Example}
\newtheorem{clm}{Claim}
\newtheorem{op}{Problem}
\newtheorem{con}{Conjecture}

\def \bp {\begin{prp} \ }
\def \ep {\end{prp}}

\def \bc {\begin{crl} \ }
\def \ec {\end{crl}}

\def \bcon {\begin{con} \ }
\def \econ {\end{con}}

\def \thm {\begin{Theorem} \ }
\def \ethm {\end{Theorem}}

\def \bl {\begin{lem} \ }
\def \el {\end{lem}}

\def \bd {\begin{defi} \ \rm }
\def \ed {\end{defi}}

\def \brm {\begin{rmk} \ }
\def \erm {\end{rmk}}

\def \bxm {\begin{xmp} \ \rm }
\def \exm {\end{xmp}}

\def \bcm {\begin{clm} \ }
\def \ecm {\end{clm}}

\def \bop {\begin{op} \ }
\def \eop {\end{op}}

\def \nmr {\begin{enumerate}}
\def \enmr {\end{enumerate}}

\def \tmz {\begin{itemize}}
\def \etmz {\end{itemize}}

\def \bsk {\bigskip}

\def \es {\emptyset}

\newcommand{\proof}{\noindent{\bf Proof.\ }}
\newcommand{\qed}{\hfill $\square$ \medskip}

\begin{document}
\title{Improved upper bounds on the domination number of graphs with minimum degree at least five}

\author{
Csilla Bujt\'as $^{a}$
\and
Sandi Klav\v zar $^{b,c,d}$
}

\date{\today}

\maketitle

\begin{center}
$^a$ Department of Computer Science and Systems Technology\\
University of Pannonia,  Veszpr\'em, Hungary \\
{\tt bujtas@dcs.uni-pannon.hu} \\

\medskip
$^b$ Faculty of Mathematics and Physics, University of Ljubljana, Slovenia\\
{\tt sandi.klavzar@fmf.uni-lj.si}

\medskip
$^c$ Faculty of Natural Sciences and Mathematics, University of Maribor, Slovenia

\medskip
$^d$ Institute of Mathematics, Physics and Mechanics, Ljubljana\\
\end{center}

\maketitle

\begin{abstract}
An algorithmic upper bound on the domination number $\gamma$ of graphs in terms of the order $n$ and the minimum degree $\delta$ is proved. It is demonstrated that the bound improves best previous bounds for any $5\le \delta \le 50$. In particular, for $\delta=5$, Xing et al.\ proved in 2006 that $\gamma \le 5n/14 < 0.3572 n$. This bound is improved to $0.3440 n$. For $\delta=6$, Clark et al.\ in 1998 established $\gamma <0.3377 n$, while Bir\'o et al.
recently improved it to $\gamma <0.3340 n$. Here the bound is further improved to $\gamma < 0.3159 n$. For $\delta=7$, the best earlier bound $0.3 088 n$ is improved to $\gamma < 0.2927 n$.
\end{abstract}


\medskip\noindent
{\bf AMS Subj. Class:}  05C69, 05C35

\newpage
\section{Introduction}


As usual, the domination number of a graph $G$ will be denoted with $\gamma(G)$ and the order of $G$ with $n$. Unless stated otherwise, the graphs considered will be connected.

A central theme in domination theory is a search for upper bounds for the domination number of graphs of given minimum degree in terms of the order of a graph. An early general bound due to Arnautov~\cite{arnautov-1974} and Payan~\cite{payan-1975} asserts that 
\begin{equation}
\label{eq:arnautov}
\gamma(G) \le \frac{n}{\delta +1}\sum_{j=1}^{\delta+1}\frac{1}{j}\,
\end{equation} 
holds for any graph $G$ of minimum degree $\delta$ and order $n$. As a consequence, 
\begin{equation}
\label{eq:alon}
\gamma(G) \le n \left( \frac{1 + \ln (\delta + 1)}{\delta +1}\right)\,.
\end{equation} 
The bound~\eqref{eq:alon} also follows from a more general result on transversals in hypergraphs due to Alon~\cite{alon-1990}. His proof is probabilistic and can also be used to infer that~\eqref{eq:alon} is asymptotically (that is, when $\delta \to \infty$) optimal. For a proof of~\eqref{eq:alon} in terms of graphs see~\cite[Theorem 10.5]{imrich-2008}.

Many investigations were done for specific values $\delta$ in order to improve the above general bounds. The first result goes back to Ore~\cite{ore-1962} who observed that if $\delta(G)\ge 1$ then $\gamma(G)\le n/2$. Blank~\cite{blank-1973}, and later independently McCuaig and Shepherd~\cite{mccuaig-1989}, followed by proving that $\gamma(G)\le 2n/5$ holds for all graphs $G$ with $\delta(G)=2$ except for seven small graphs ($C_4$, and six graphs on seven vertices). For graphs $G$ with $\delta(G)=3$, Reed~\cite{reed-1996} proved his celebrated result: $\gamma(G) \le 3n/8$. This bound is sharp as there exist cubic graphs of order 8 with domination number 3. For cubic graphs Reed's result was further improved by Kostochka and Stodolsky~\cite{kostochka-2009} by proving that as soon as a cubic graph has at least nine vertices, $\gamma(G)\le 4n/11$ holds. For additional closely related interesting results see~\cite{kral-2012,lo-2008,lo-2012}. We also add that in~\cite{volk-2006} an upper bound of different nature on the domination number in terms of order and minimum degree is given. 

The above bounds of Ore, Blank, and Reed have the same shape: 
\begin{equation}
\label{eq:simple}
\gamma(G) \le \frac{n\delta}{3\delta-1}\,,
\end{equation}
where $1\le \delta \le 3$ is the minimum degree of $G$ and in the case $\delta=2$ we skip the seven exceptions. Actually,~\eqref{eq:simple} holds for any minimum degree $\delta\ge 1$ as conjectured in~\cite[p.~48]{haynes-1998}. For $\delta\ge 6$, the bound~\eqref{eq:arnautov} is better than~\eqref{eq:simple}, while for $\delta = 4, 5$ the bound~\eqref{eq:simple} was proved in~\cite{sohn-2009,XSC}, respectively.

Clark et al.~\cite{CSSF} proved the following stronger result 
\begin{equation}
\label{eq:clark}
\gamma(G) \le n\left( 1 - \prod_{j=1}^{\delta +1} \frac{j\delta}{j\delta + 1} \right)\,,
\end{equation}
which is better than~\eqref{eq:simple} for any $\delta\ge 6$ and better than~\eqref{eq:arnautov} for any $\delta\ge 5$. 

Recently, Bir\'o et al.~\cite{BCDS} further improved the bound~\eqref{eq:clark} by proving that 
\begin{equation}
\label{eq:biro}
\gamma(G) \le n\left( 1 - \frac{\delta^2-\delta+1}{1 + \delta \prod_{j=1}^{\delta - 1} \left( 1 + \frac{\delta+1}{j\delta}\right)}\right)\,.
\end{equation}
The bound~\eqref{eq:biro} is better than~\eqref{eq:clark} for any $\delta$ and is better than~\eqref{eq:arnautov} for $\delta\ge 6$. The present state of the art up to the minimum degree 7 is summarized in Table~\ref{table:best-up-to-now}.

\begin{table}[ht!]
\begin{center}
\begin{tabular}{|c||c|c|c|c|c|c|c|}
\hline
$\phantom{\Big\rvert}\delta$ & 1 & 2 & 3 & 4 & 5 & 6 & 7 \\
\hline
$\phantom{\Bigg\rvert}\gamma(G)\le$ &
$\displaystyle{\frac{n}{2}}$ &
$\displaystyle{\frac{2n}{5}}$ &
$\displaystyle{\frac{3n}{8}}$ &
$\displaystyle{\frac{4n}{11}} < 0.3637n$ &
$\displaystyle{\frac{5n}{14}}<0.3572n$ &
$0.3340n$ &
$0.3089n$ \\
\hline
\end{tabular}
\end{center}
\label{table:best-up-to-now}
\caption{Best present upper bounds, where for $\delta = 2$ seven exceptional graphs are not included}
\end{table}

The bounds from Table~\ref{table:best-up-to-now} are sharp for $\delta \le 3$, while the sharpness for $\delta = 4$ is not known. In this paper we improve the best earlier bounds for any minimum degree $\delta$, $5\le \delta\le 50$. 

In the next section we present our main theorem and explain how the new bounds can be derived from it. In the subsequent section the main result is proved. Its proof idea is a variation of the recent approach from the theory of the domination game~\cite{bujtas-2013,bujtas-2014}.

\section{Main result and its consequences}

Before stating our main theorem, we emphasize that the result is quite technical. We hence ask the reader to judge its usefulness by the results presented after its statement. 

\thm
\label{thm:main}
Let $G$ be a graph of minimum degree $d\ge 5$. Let $a$, $s$, and $b_1,\ldots, b_d$ be positive numbers such that
\begin{enumerate}
\item[(i)]  $0 \le b_d-b_{d-1} \le b_{d-1}-b_{d-2}\le \dots\le b_2-b_1 \le b_1$, and $b_d\le a$,
\item[(ii)]  $(d+2)a-(d+1)b_d \ge s$,
\item[(iii)] $(d+1)a-db_{d-1} \ge s$,
\item[(iv)] for $2\le i\le d-1$,

$(d-i+2)a+(d-i+2)(i-1)b_{d-i+2}-(d-i+2)(i-1)b_{d-i+1}-(d-i+1)b_{d-i} \ge s$,
\item[(v)] $2a+2(d-1)b_2-2(d-1)b_1\ge s$, and
\item[(vi)] $a+d b_1\ge s$.
\end{enumerate}
Then
$$\gamma(G)\le \frac{a}{s}\; n.$$
\ethm

To apply Theorem~\ref{thm:main}, we need to explicitly determine the values $a$ and $s$. For this sake, we first solve the system of $d+1$ linear equations which correspond to conditions (iii)-(vi) of the theorem, where we set equality instead of the inequality, and select a fixed $s$, say $s=1$. The obtained system always has a unique solution yielding the values of $a$, $b_1, \ldots, b_d$. Indeed, first $b_1$ can be computed as a function of $a$ from (vi), then $b_2, \ldots ,b_d$ can be computed as functions of $a$ from (iv) and (v), finally $a$ is determined by (iii).  After $a$, $b_1, \ldots, b_d$ are computed for a fixed $d$, we also check that (i) and (ii).

For graphs with $\delta(G) \in \{5, 6, 7\}$  computations give: 

\bc
(i) If $\delta(G) = 5$, then $$\displaystyle{\gamma(G) \le \frac{2671}{7766}\,n}\,.$$

(ii) If $\delta(G) = 6$, then $$\displaystyle{\gamma(G) \le \frac{1702}{5389}\,n}\,.$$

(iii) If $\delta(G) = 7$, then $$\displaystyle{\gamma(G) \le \frac{389701}{1331502}\,n}\,.$$
\ec

\proof
As described above, we have performed the corresponding computations and obtained the following values, where at the end we have changed $s$ from 1 to an appropriate integer:
\begin{itemize}
\item $\delta(G) = 5$: $a=2671$,    $b_5=1751$,    $b_4=1652$,    $b_3=1521$,    $b_2=1322$,    $b_1=1019$, $s=7766$; 
\item $\delta(G) = 6$: $a=1702$,    $b_6=1137$,    $b_5=1087.5$,  $b_4=1024$,    $b_3=939$, $b_2= 813$, $b_1= 614.5$, $s=5389$;  
\item $\delta(G) = 7$: $a=1169103$,  $b_7=793539$,  $b_6=765474$,  $b_5=730945$,   $b_4=686892$,  $b_3=627951$,  $b_2=541654$,  $b_1=403629$,
  $s=3994506$.
\end{itemize}
The reader can check that in all the cases the conditions of Theorem~\ref{thm:main} are fulfilled. 
\qed

We add that Theorem~\ref{thm:main} also holds for graphs of minimum degree $3$ and $4$, but for these cases we do not obtain better bounds on $\gamma(G)$ than the best earlier ones. In Table~\ref{table:summarized} the upper bound of Theorem~\ref{thm:main} is compared with the Arnautov bound~\eqref{eq:arnautov}, the upper bound~\eqref{eq:simple}, and the Bir\'o et al.\ bound~\eqref{eq:biro}. For the sake of compactness we list only values up to $\delta=20$,  but we have computed all the values up to $\delta = 50$ and found that Theorem~\ref{thm:main} leads to best bounds up to date. 

\begin{table}[ht!]
\begin{center}
\begin{tabular}{|c||c|c|c|c|}
\hline
$\phantom{\Big\rvert}\delta$ & \eqref{eq:arnautov} & \eqref{eq:simple} & \eqref{eq:biro} & Theorem~\ref{thm:main} \\
\hline
 5 & 0.380556 & 0.357143 & 0.364253 & 0.343935 \\
\hline
 6 & 0.350000 & 0.352941 & 0.333938 & 0.315829 \\
\hline
 7 & 0.324107 & 0.350000 & 0.308805 & 0.292678 \\
\hline
 8 & 0.301984 & 0.347826 & 0.287619 & 0.273213 \\
\hline
 9 & 0.282897 & 0.346154 & 0.269496 & 0.256566 \\
\hline
10 & 0.266270 & 0.344828 & 0.253796 & 0.242128 \\
\hline
11 & 0.251656 & 0.343750 & 0.240046 & 0.229463 \\
\hline
12 & 0.238709 & 0.342857 & 0.227891 & 0.218244 \\
\hline
13 & 0.227152 & 0.342105 & 0.217057 & 0.208223 \\
\hline
14 & 0.216771 & 0.341463 & 0.207331 & 0.199207 \\
\hline
15 & 0.207389 & 0.340909 & 0.198545 & 0.191045 \\
\hline
16 & 0.198866 & 0.340426 & 0.190562 & 0.183614 \\
\hline
17 & 0.191086 & 0.340000 & 0.183273 & 0.176815 \\
\hline
18 & 0.183953 & 0.339623 & 0.176588 & 0.170566 \\
\hline
19 & 0.177387 & 0.339286 & 0.170430 & 0.164801 \\
\hline
20 & 0.171321 & 0.338983 & 0.164738 & 0.159462 \\
\hline
\end{tabular}
\end{center}
\label{table:summarized}
\caption{Comparison of Theorem~\ref{thm:main} with earlier upper bounds}
\end{table}

\section{Proof of Theorem~\ref{thm:main}}

Before actual proof we explain its main idea and the source for it. The domination game introduced in~\cite{bresar-2010} is played by two players who alternate in choosing vertices of a graph such that each chosen vertex enlarges the set of vertices dominated so far. The aim of one player is that the graph is dominated in as few steps as possible, while the aim of the other player is just the opposite. Kinnersley, West, and Zamani~\cite{bill-2013} conjectured that no domination game lasts more that $3/5$ of the order of the graph. To attack this conjecture a method was introduced in~\cite{bujtas-2013,bujtas-2014} that, roughly, aims to weight the vertices of a graph and accordingly change the values during the course of the game. In the present situation we have only one player, and the idea can be modelled as a greedy algorithm for estimating the domination number of a graph in which the greedy criteria is designed according to the weights assigned to the vertices. We note in passing that a complementary approach was taken in~\cite{bresar-2014}, more precisely, the authors considered legal dominating sequences of maximum length (and named the length of it Grundy domination number) which can be understood as the domination game with one player but now the player wants the game to last as long as possible. 

In the rest of the section we prove Theorem~\ref{thm:main}. In the proof, we consider a graph $G=(V,E)$ of minimum degree $d\ge 5$, and construct a  dominating set  selecting its vertices one-by-one.  At the beginning, we set $D=\es$ and then, in each step the vertex chosen is put into $D$. It will be clear that the procedure analyzed here can be interpreted as a greedy algorithm; that is, in each step we select a vertex which dominates the most vertices undominated before.

 Our main tool
 in the proof  is a value assignment $p: \; V \rightarrow
 \mathbb{R}$ which always relates to the current set $D$; hence, in each
 step of the algorithm we have some vertices $v$ whose value $p(v)$
 changes.
 We also make  distinction between white, blue and red
 vertices, and for each vertex  $v$ its white degree $\deg_W(v)$ ({\em W-degree}) is
 just the number of its white neighbors. The blue degree $\deg_B(v)$
 is defined analogously.
  We use the following terminology and notations:
 \tmz
 \item A vertex  is {\it white} if it is undominated by $D$.
 Every white vertex is assigned with the same value denoted by $a$.
 \item A vertex $v$ is {\it blue} if it is dominated but has at
 least one white neighbor. The value assigned to $v$ depends on the
 number $\deg_W(v)$ of its white neighbors as follows: 
 \tmz
 \item if $\deg_W(v) \ge d$, then $p(v)=b_d$;
 \item if $\deg_W(v) = i$ for an $1\le i \le d-1$, then $p(v)=b_i$.
 \etmz
 \item A vertex $v$ is {\it red}
 and $p(v)=0$ if each vertex from $N[v]$ is dominated by $D$.
 \etmz

 By definitions, no white vertex $v$ has a red neighbor, hence
 $\deg_B(v)\ge d-\deg_W(v)$. Moreover, for any vertex $u$, the
 number  of its white neighbors cannot increase during the
 procedure. Especially, if $\deg_W(v)\le k$ holds for the vertex
 $v$ at a moment, this remains true in each later step  of the
 algorithm. Note also that a blue vertex always has at least one white neighbor.
 \bsk

 In every step, the {\it value of the graph $G$} is just the sum
 $p(G)=\sum_{v \in V} p(v)$. By the way $p$ is defined and by the conditions of the theorem   given in part $(i)$ which imply that $a\ge b_d\ge b_{d-1}\ge \cdots \ge b_2\ge b_1\ge 0$, it follows that $p(G)$ decreases in each step. This reduction
 is called the {\it gain} of the step. Clearly, $p(G)=na$
 when the algorithm starts  and $p(G)=0$ at the end when $D$ is a dominating set. Thus, once
 we prove that the gain is at least $s$ in each step, the desired
 inequality $\gamma \le an/s$ will follow.

 The process of constructing the dominating set is divided into
 $d+2$ phases some of which might be empty.

 \tmz
  \item If there exists a vertex   whose closed neighborhood
  contains at least $d+2$ white vertices, then the next choice
  belongs to Phase 0.
  \item If there exists a vertex $v$  whose closed neighborhood
  contains  $d-i+2$ white vertices (for $1 \le i \le d+1$), but no other vertex $w$ has more than
  $d-i+2$ white vertices in $N[w]$, then the next choice
  belongs to Phase $i$.
   \etmz

Clearly, each choice belongs to exactly one phase. Further, in Phase
$i$ every white vertex $v$ has $\deg_W(v) \le d-i+1$ and every blue
vertex $u$ has $\deg_W(u) \le d-i+2$. The construction of the
dominating set can be done in two different ways:
  \tmz
    \item[$(1)$] For every $0 \le i\le d+1$, in each step of Phase $i$ we choose a vertex whose closed
  neighborhood contains at least $d-i+2$ white vertices. Under this
  condition we prefer to select a white vertex.
  \item[$(2)$] In each step, a vertex with the possible largest gain is chosen
  \etmz
  For the sake of simplicity, we will consider  strategy $(1)$, but note
  that our proof remains  valid  when $(2)$ is followed.

  \bsk

 Now, we are ready to prove that following the greedy strategy $(1)$, the gain of each step is at least $s$.

   \paragraph{Phase 0.} If there is a white vertex $v$ with at least $d+1$ white
  neighbors, then selecting $v$, its color turns from white to red
  and its value $p(v)$ decreases by $a$. Moreover, each white
  neighbor $u$ of $v$ turns from white to either blue or red,
  which means a decrease of at least $a-b_d$ in $p(u)$. Thus, the
  choice of $v$ reduces $p(G)$ by at least
  $(d+2)a-(d+1)b_d$. By $(ii)$, this gain is not smaller than $s$.
  The situation is similar if we have a blue vertex $v$ with
  $\deg_W(v)\ge d+2$. Selecting $v$, $p(v)$ decreases by exactly $b_d$
  and the value of each white neighbor is reduced by at least
  $(a-b_d)$.

   \paragraph{Phase 1.} First,  assume there exists a white vertex $v$
    with $\deg_W(v)=d$. In Phase 1, each white vertex has at most
    $d$ white neighbors. Hence, after the choice of $v$, the value of each (originally) white
    neighbor decreases by at least $a-b_{d-1}$. Then, the gain is at least $s$ by the  inequality $(iii)$.
    If every white vertex has at most $d-1$ white neighbors, but the
    choice belongs to Phase 1, then we have a blue vertex $v$ with
    $d+1$ undominated neighbors. When such a blue vertex $v$ is selected,
    the new value of its any white neighbor  cannot
    exceed $b_{d-1}$. Thus, the gain is at least $b_d+(d+1)(a- b_{d-1}) \ge
    (d+1)a-db_{d-1} \ge s$.

    \bsk

    \paragraph{Phase $\bm{i}$, $\bm{2 \le i \le d+1}$.} In Phase $i$ each  
    white vertex is of W-degree at most $d-i+1$ and each blue vertex is of 
    W-degree at most $d-i+2$. This fact together with our condition $(i)$ implies
     that whenever the W-degree of a
    blue vertex is decreased by $\ell$, its value decreases by at
    least $\ell(b_{d-i+2}-b_{d-i+1})$ (for $2\le i\le d$).

    Then, if $2 \le i \le d-1$ and a white vertex $v$
    with white neighbors $u_1,\dots,u_{d-i+1}$ is chosen, the
    following changes occur in the values of vertices:
    \tmz
    \item $p(v)$ is reduced by $a$;
    \item for every $1\le j\le d-i+1$, $p(u_j)$ is reduced by at
    least $a-b_{d-i}$;
    \item as each of $v,u_1,\dots,u_{d-i+1}$ has at least $i-1$ blue
    neighbors, the decrease in the sum of the values of blue
    neighbors is not smaller than
    $(d-i+2)(i-1)(b_{d-i+2}-b_{d-i+1})$.
    \etmz
    Then, condition $(iv)$ of the theorem  ensures that the total
     decrease is at least $s$.
     \bsk

     If there exists no white vertex with W-degree at least $d-i+1$, but
     we have a blue vertex $v$ with $d-i+2$ white neighbors
     $u_1,\dots,u_{d-i+2}$, then choosing $v$ we have the following
     changes:
     \tmz
      \item $p(v)$ decreases by $b_{d-i+2}$;
    \item for every $1\le j\le d-i+2$, $p(u_j)$ is reduced by at
    least $a-b_{d-i}$;
    \item as each of $u_1,\dots,u_{d-i+2}$ has at least $i-1$ blue
    neighbors different from $v$, we have  additional decrease of at
    least
    $(d-i+2)(i-1)(b_{d-i+2}-b_{d-i+1})$.
     \etmz
     That is, the total change in $p(G)$ is at least
     $$b_{d-i+2}+(d-i+2)(a-b_{d-i})+
     (d-i+2)(i-1)(b_{d-i+2}-b_{d-i+1}). $$
    By $(iv)$ it is at least $s$, as also $b_{d-i+2}\ge b_{d-i}$
    must hold due to  $(i)$.

     \bsk

     Hence, the gain is at least $s$ in every step which belongs to Phase
     $i$ where $2 \le i \le d-1$, and
     by a similar argumentation we obtain that the same is true
     for Phase $d$ by the condition $(v)$ and for Phase $(d+1)$ by the condition $(vi)$.
     \bsk

     It follows that a dominating set is obtained in at most $an/s$ steps
     and consequently 
     $$\gamma(G)\le \frac{a}{s}\; n.$$

\section{Concluding remarks}

It would be interesting to see how good are the upper bounds obtained in this paper. In this respect we recall that for given $n$ and $\delta$, Clark and Dunning~\cite{clark-1997} defined $\gamma(n,\delta)$ to be the maximum domination number of an arbitrary graph (that is, connected or disconnected) of order $n$ and minimum degree $\delta$. They determined $\gamma(n,\delta)$ for all $n$ and $1\le \delta \le 3$ (where the examples constructed are mainly disconnected) and most of the values for $n\le 14$. From our point of view it would be most interesting to find connected graphs of an arbitrary order that are close to the obtained upper bounds. 

To conclude the paper we add that we strongly believe that Theorem~\ref{thm:main} improves the present best upper bounds for any $\delta \ge 5$.

\section*{Acknowledgements}
Research of the first author was supported by the European Union and
Hungary through the projects T\'AMOP-4.2.2.C-11/1/KONV-2012-0004 and the Campus Hungary B2/4H/12640. The second author was supported by the Ministry of Science of Slovenia under the grant P1-0297.

%
\end{document}